\documentclass[11pt]{article}
\usepackage{mathpazo}

\usepackage{amssymb,amsmath,amsthm,amsfonts,amscd,latexsym, dsfont, color, xcolor}
\usepackage{graphics}
\usepackage{hyperref} 
\usepackage{array} 
\usepackage[pdftex]{graphicx}
\hypersetup{colorlinks}
\usepackage{wrapfig}

\usepackage{fancyhdr}

\definecolor{indigo}{RGB}{51,0,102}
\definecolor{brightpurple}{RGB}{102,0,153}
\definecolor{fuchsia}{RGB}{180,51,180}
\definecolor{jolightpurple}{RGB}{188,171,240}

\hypersetup{colorlinks,
linkcolor=brightpurple,
filecolor=brightpurple,
urlcolor=indigo,
citecolor=fuchsia}

\usepackage[margin=1in]{geometry}
\normalbaselines
\newcommand{\Z}{\mathbb{Z}}
\newcommand{\R}{\mathbb{R}}

\begin{document}
\vspace{-10pt}
\section{Introduction}
 \emph{Symplectic and contact topology} is an active area of mathematics that combines ideas from dynamical systems, analysis, topology, several complex variables, as well as differential and algebraic geometry.  Symplectic and contact structures first arose in the study of \emph{classical mechanical systems}, allowing one to describe the time evolution of both simple and complex systems such as springs, planetary motion and wave propagation \cite{Ar}.   Understanding the evolution and distinguishing transformations of these systems led to the development of global invariants of symplectic and contact manifolds.  

 The equations of motion in classical mechanics are determined by the notion of a conserved quantity, \emph{energy}.  A related quantity is \emph{action}, which is minimized by solutions to the equations of motion.  For a closed system, such as the Kepler problem whose solutions describe paths of planets orbiting the sun,  the energy is the sum of the kinetic and potential energy in the system, and the action is given by the (minimized) mean value of kinetic minus potential energy.  Symplectic and contact structures emerge as we investigate these systems by unpacking the information hidden in the notions of energy and action.

The position of a particle in a mechanical system is a point $x = (x_1,..,x_n)$ in Euclidean space, and the vector space $\R^n$ defined by these coordinates is called the \emph{configuration space}.  The position and momentum of a particle allows us to predict the particle's motion at all future times within a system.  The \emph{phase space} of a system is precisely this space that represents all possible states of the system, consisting of both the position and momentum of a particle.  In the case that there are $n$ degrees of freedom, the phase space is $\R^{2n}$.     The assumption that the trajectories of a particle $x(t)$ minimize an action functional gives rise to a system of $n$ second-order differential equations called the \emph{Euler-Lagrange equations}, discovered in 1808 by Joseph-Louis Lagrange \cite{L}.  

These equations grew out of Lagrange's observation that the possible elliptic motions of a single planet under the sun's gravitational pull can be described by six real parameters.   However, the influence of other planets perturbs this ellipticity.  In order to describe the variation, one must study the derivatives of these real parameters.  These three equations are extremely complicated, but they can be simplified by introducing \emph{Lagrange brackets}, which are combinations of the derivatives with respect to position and velocity at fixed time.   Lagrange then showed that these equations can be transformed into what is now known as a \emph{Hamiltonian system} of six first-order differential equations that conserve energy \cite{L2}.  At the time, his notion of energy was a ``disturbing function," which described the variance from elliptic motion.  Moreover, these Lagrange brackets turn out to be the coefficients of humanity's oldest symplectic structure \cite{wlectures}.   
 
 \begin{wrapfigure}{l}{0.47\textwidth}
    \vspace{-5pt}
    \includegraphics[width=0.47\textwidth]{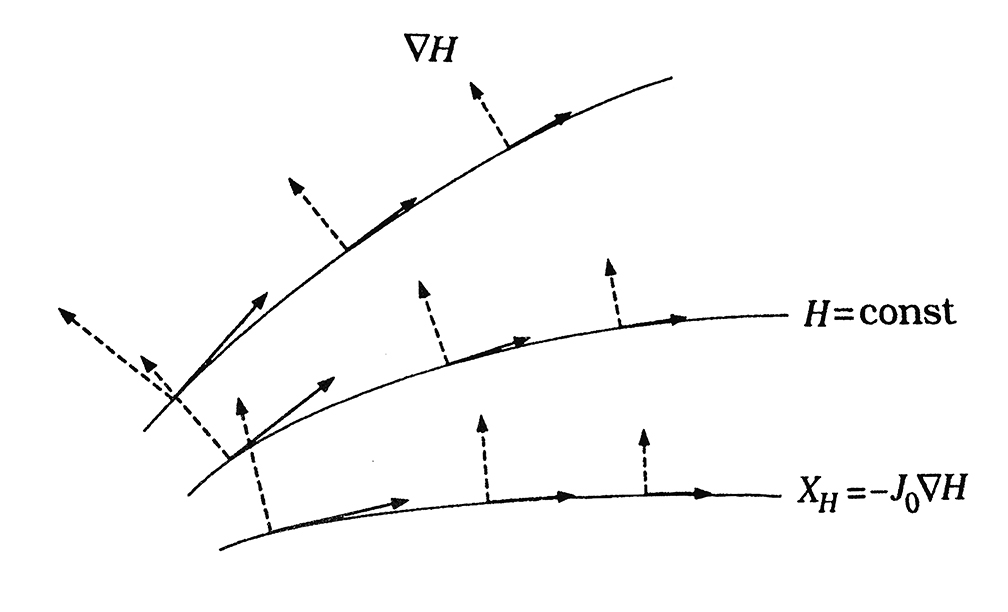}
        \vspace{-30pt}
  \caption{ {\footnotesize The symplectic gradient $X_H$, \cite{MS1} }}
  \label{ham}
  \vspace{-10pt}
\end{wrapfigure} 

In the mid 1800s, William Rowan Hamilton and Carl Jacobi realized the theoretical consequences of Lagrange's work, in particular that the $n$ Euler-Lagrange equations can be transformed into a Hamiltonian system of $2n$ equations \cite{W}.  The Hamiltonian system is governed by the conservation of an energy function, called the Hamiltonian function $H(x,y)$, which defines the \emph{Hamiltonian vector field} $X_H$.  The flow lines of this vector field are solutions to Hamilton's equations of motion,
  \[
 \dot{x} = \dfrac{\partial H}{\partial y}, \ \ \ \dot{y} = -\dfrac{\partial H}{\partial x}.
 \] 
 
In the coordinates $z=(x_1,...x_n,y_1,...y_n) \in \R^{2n}$ the Hamiltonian system can be written in the form of a system of $2n$ differential equations, 
\[
J_0 \dot{z} = \nabla H(z).
\]
where  $\nabla H$ denotes the gradient of $H$ and $J_0$ is the $2n \times 2n$ matrix
\[
J_0=\left( \begin{array}{cc} 0 &- \mathds{1} \\ \mathds{1} &0  \\\end{array} \right).
\]
The Hamiltonian vector field or \emph{symplectic gradient} of $H$, seen in Figure \ref{ham}, is defined by
\[
X_{H} = -J_0 \nabla H : \R^{2n} \to \R^{2n}.
\]

Systems whose Hamiltonian function explicitly depends on time, such as those describing the motion of a charged particle in a time-dependent electric field, use \emph{extended phase space}, which includes the $2n$-phase space plus the time variable.  Extended phase space results in the notion of a \emph{contact structure}.  In this setting, solutions to equations of motion yield flows of a Hamiltonian-like vector field, called the Reeb vector field.  

Contact structures appear naturally in other areas of mathematics and physics, including thermodynamics \cite{atherm}.  In particular, contact geometry allows one to understand geodesic flow on the tangent bundle of a Riemannian manifold.  Geodesics are locally the shortest distance between points, where distance is defined in terms of a metric intrinsic to a manifold.  An $n$-dimensional \emph{manifold} is a smooth object that locally looks like $\R^n$.\footnote{For example, the surface of a donut or beach ball is a 2-manifold.  If we cut out a small piece of either surface and ``zoomed in'' it would look like a flat sheet of paper, e.g. $\R^2$.}   One can interpret a Riemannian manifold as a model for an optical medium, in which case geodesics with respect to the metric correspond to light rays.  This in turn yields \emph{Huygens's principle}, which states that every point on a wavefront is a source of wavelets, which spread forward at the same speed.

\section{Symplectic and contact manifolds}
 
 \begin{wrapfigure}{r}{0.26\textwidth}
  \vspace{-65pt}
 \includegraphics[width=0.25\textwidth]{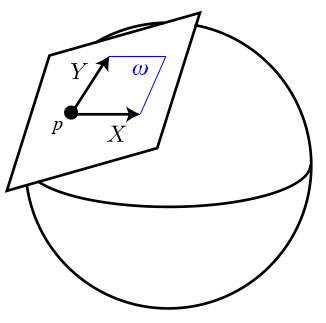}
  \vspace{-10pt}
  \caption{\footnotesize At the infinitesimal level, $\omega$ measures oriented area spanned by vectors $X$ and $Y$ at a point $p$.}
\vspace{-10pt}
  \label{omega}
\end{wrapfigure}
 To study more general even dimensional Hamiltonian systems we need to allow symplectic manifolds to serve as the phase space.     In classical mechanics, replacing the standard $2n$-dimensional phase space with a $2n$-manifold results in a canonical \emph{symplectic structure} on the manifold, reflecting the conservation of energy.   Formally, a symplectic form $\omega$ is a closed nondegenerate 2-form.   It allows one to measure two dimensional area in a well-defined way, as seen in Figure \ref{omega}, and as a result forces symplectic manifolds to be even dimensional. Using the symplectic form one can define the \emph{Hamiltonian vector field}, $X_H$, on a symplectic manifold by
 \[
 \omega(X_H, \cdot) = dH(\cdot).
 \]

The name symplectic arose in 1939 due to Hermann Weyl, who studied the symplectic linear group.  This group manifests itself when one studies the \emph{canonical transformations}\footnote{Now such isomorphisms are called symplectomorphisms, due to Souriau's contributions \cite{SR}.} of a Hamiltonian system, which are changes of coordinates that preserve  Hamilton's equations.  Weyl recalls in a footnote on page 165 \cite{weyl}, ``The name \emph{complex group} formerly advocated by me in allusion to line complexes, as these are defined by the vanishing of antisymmetric bilinear forms, has become more and more embarrassing through collision with the word \emph{complex}, [a Latin adjective],  in the connotation of complex number. I therefore propose to replace it by the corresponding Greek adjective \emph{symplectic}."  
 
\begin{wrapfigure}{r}{0.55\textwidth}
  \vspace{-15pt}
    \includegraphics[width=0.55\textwidth]{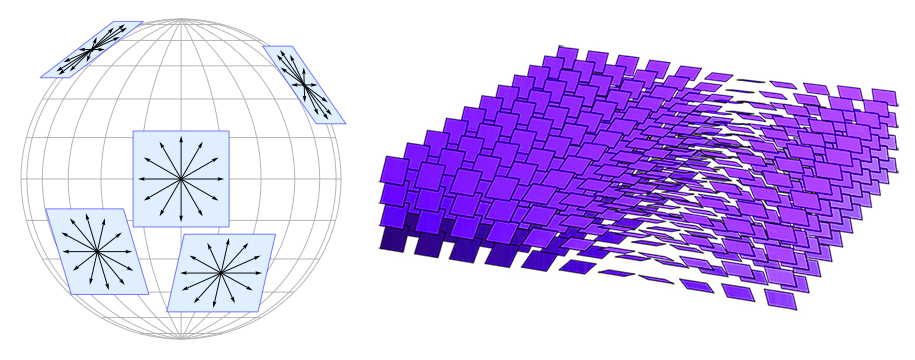}
  \vspace{-25pt}
  \caption{ \mbox{\footnotesize An integrable (right) \& contact (left) structure on $\R^3$.}}
  \vspace{-15pt}
  \label{fig3}
\end{wrapfigure} 

 Many contact manifolds arise as hypersurfaces or boundaries of symplectic manifolds, and the geometry of contact and symplectic manifolds is closely intertwined. A \emph{contact structure} $\xi$ is a maximally nonintegrable hyperplane distribution. In three dimensions, this means that the planes of $\xi$ twist so much that even locally there is never a surface whose tangent planes are all contained in $\xi$, which is in contrast to the notion of an integrable hyperplane distribution, seen in Figure \ref{fig3}.  An integrable hyperplane distribution is one in which all the planes are given by tangent planes of a submanifold.   Any 1-form $\alpha$ whose kernel defines a contact structure is called a \emph{contact form}.

   \begin{wrapfigure}{L}{.45\textwidth}
  \vspace{-20pt}
  \begin{center}
    \includegraphics[width=0.45\textwidth]{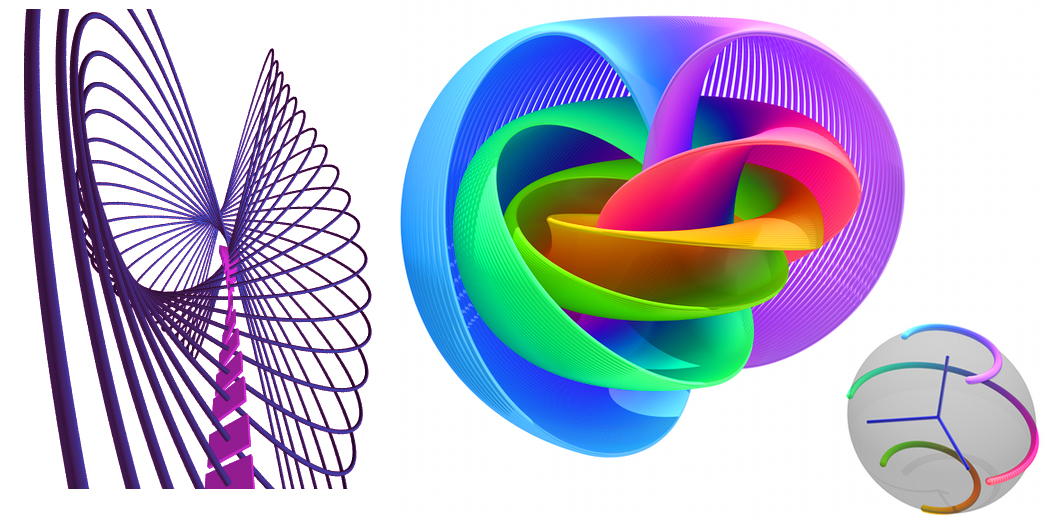}
  \end{center}
  \vspace{-10pt}
  \caption{\footnotesize The flows of two Reeb vector fields; the right is on $S^3$ and is parametrized by $S^2.$ Credit: Patrick Massot (left), Niles Johnson (right)}
  \label{fig:reeb}
\end{wrapfigure}  

 The \emph{Reeb vector field} $R_\alpha$ depends on the choice of contact form $\alpha$ and is defined by 
\[
\alpha(R_\alpha)=1, \ \ \ d\alpha(R_\alpha, \cdot) = 0.
\]
The flow of  $R_\alpha$ preserves the form $\alpha$ and hence the contact structure $\xi$. It can also follow very complex patterns, as in Figure \ref{fig:reeb}. 

Moreover, the flows Reeb vector fields of different contact forms defining the same contact structure may have wildly different properties. 

An interesting result about symplectic and contact manifolds is Darboux's theorem, which states that locally all contact structures look like the kernel of the standard contact form on $\R^{2n+1},$
\[
\xi_0 = \mbox{ker}\ \alpha_0 = \mbox{ker} \left( dz +\sum_{i=1}^n x_idy_i \right),
\]
and that locally all symplectic forms look like the standard symplectic form on $\R^{2n},$
\[
\omega_0 = \sum_{i=1}^n dx_idy_i.
\]
Hence, there can be no local invariants of symplectic and contact manifolds, a stark contrast to Riemannian geometry where the notion of curvature provides local invariants. In the symplectic realm, the absence of local invariants means that there is an  infinite dimensional group of diffeomorphisms that preserve the symplectic structure and a discrete set of nonequivalent \emph{global} symplectic structures in each cohomology class.  Analogously in the contact realm, there is an infinite dimensional group of diffeomorphisms that preserve the contact structure and a discrete set of nonequivalent global contact structures in each planar homotopy class. 

The ability to distinguish contact structures in a planar homotopy class is not obvious.  One of the first results along these lines is the celebrated theorem of Yakov Eliashberg from 1989 \cite{Eot}, which states that the 3-sphere admits two homotopy classes of contact structures which are homotopic as plane fields but which are not homotopic via contact structures. 

\begin{wrapfigure}{R}{.505\textwidth}
  \vspace{-15pt}
  \begin{center}
    \includegraphics[width=0.5\textwidth]{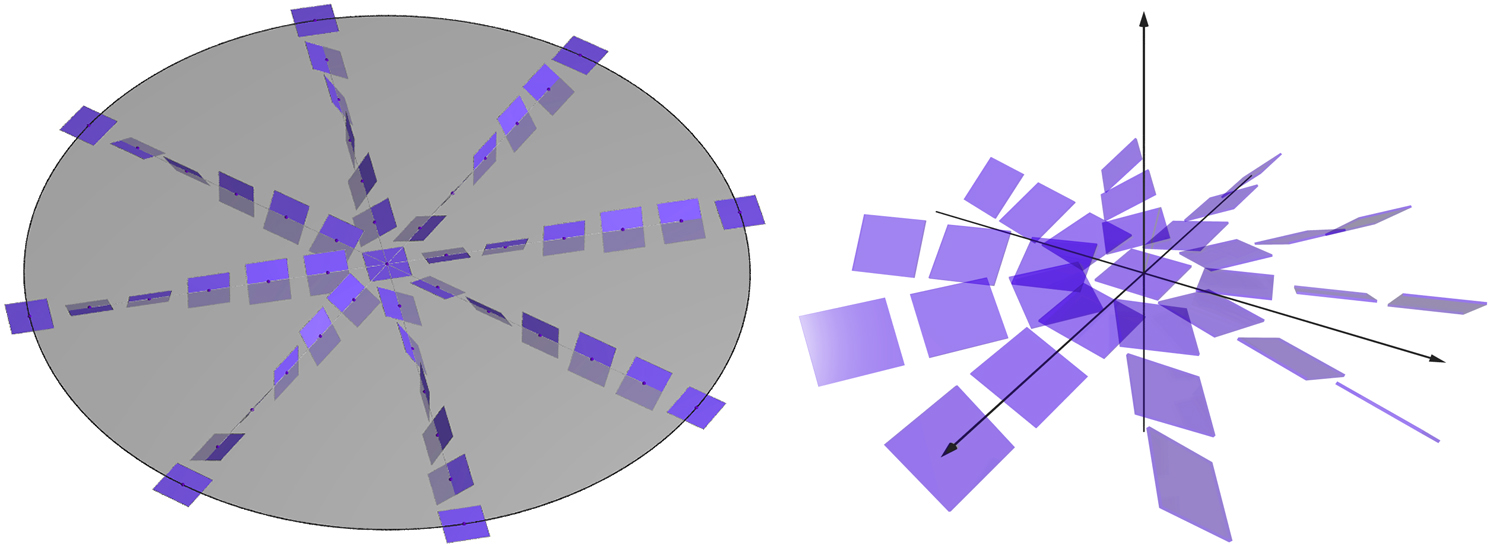}
  \end{center}
  \vspace{-20pt}
  \caption{\footnotesize The overtwisted (left, Patrick Massot) and standard (right, Otto van Koert) contact structure.}
    \vspace{-5pt}
  \label{OT}
\end{wrapfigure} 

\noindent One of these structures is the \emph{standard structure}, given in cylindrical coordinates $(r,\theta,z) \in \R^3$ by 
\[
\xi_{std} = \mbox{ker }\alpha_0  = \mbox{ker }(dz + r^2 d\theta),
\]
and the other is the \emph{overtwisted contact structure},
\[
\xi_{OT} = \mbox{ker } (\cos r dz + r \sin r d\theta ).
\]
These are visualized in the $z=0$-plane in Figure \ref{OT}.  Both $\xi_{std}$ and $\xi_{OT}$ are horizontal along the $z$-axis and along any ray they both turn counterclockwise as one moves outward from the $z$-axis.  However, the rotation angle of $\xi_{std}$ approaches (but never reaches) $\pi/2$, while the contact planes of $\xi_{OT}$ make infinitely many complete turns. 

\section{From Rabinowitz to Floer: the evolution of variational methods}
A closed orbit of a vector field $X$ on a manifold $M$ is a map,
\[
\gamma: \R /T\Z \to M,
\]
for some $T>0$, which satisfies the ordinary differential equation
\[
\dot{\gamma}(t) = X(\gamma(t)).
\]
One is then led to wonder when a (smooth) vector field $X$ on a closed manifold $M$ admits a closed orbit.  For some special three manifolds like the 3-torus, it is easy to construct vector fields with no closed orbit.  On the other hand, when $M$ is the 3-sphere, this question turns out to be incredibly difficult and not always possible; see \cite{Hu} for a brief history.  

The \emph{Weinstein conjecture} is one of the most famous questions in regard to the existence of periodic orbits \cite{weinsteinconj}.  It originated from work in the 1970s by Alan Weinstein, who demonstrated the existence of periodic orbits on convex compact hypersurfaces in $\R^{2n}$ \cite{weinstein}, and Paul Rabinowitz, who demonstrated the existence of periodic orbits on star shaped hypersurfaces in $\R^{2n}$ \cite{rabinowitz1}-\cite{rabinowitz3}.  In reading Rabinowitz's papers, Weinstein realized that there was a simple geometric feature common in the different results, namely that we he called \emph{contact type}, which is a special contact hypersurface in a symplectic manifold.   Weinstein's realization connected the existence of periodic orbits of Hamiltonian systems to contact geometry, spurring further interest in the study of contact manifolds. 
\medskip

\textbf{The Weinstein conjecture:}  Let $(M, \xi)$ be a closed co-oriented contact manifold. Then for any contact form $\alpha$ for $\xi$, the Reeb vector field $R_\alpha$ admits a closed periodic orbit.

\medskip

At the same time, Rabinowitz's paper \cite{rabinowitz2} had a profound effect on a young graduate student, Helmut Hofer.  Helmut reminisced at his sixtieth birthday conference:
\begin{quotation}
Why did I come into symplectic geometry?  I had the flu, and the only thing to read was a copy of Rabinowitz's paper where he proves the existence of periodic orbits on star-shaped energy surfaces \cite{rabinowitz2}. It turned out to contain a fundamental new idea, which was to study a different action functional for loops in the phase space rather than for Lagrangians in the configuration space.\footnote{In 1976, Moser wrote that this action functional was ``certainly not suitable for an existence proof, \cite[(1.5)]{moser}."  Rabinowitz, on the other hand, showed more optimism than his former advisor in 1977, \cite[Remark 4.44]{rabinowitz1}.}  Which actually if we look back, led to the variational approach in symplectic and contact topology, which is reincarnated in infinite dimensions in Floer theory and has appeared in every other subsequent approach. The flu turned out to be really good.    
\end{quotation}

\noindent This variational approach led to further progress by Claude Viterbo in 1987 for hypersurfaces of contact type in $\R^{2n}$ \cite{viterbo}, which was extended further by Hofer-Viterbo (\cite{HV} in 1987), Hofer-Zehnder (\cite{HZ} in 1988), and Struwe (\cite{S} in 1990).

Meanwhile, the Arnold conjecture haunted the dreams of geometers.

\medskip
\indent \textbf{The Arnold conjecture.} A symplectomorphism on a closed symplectic manifold that is generated by a time-dependent Hamiltonian vector field should have at least as many fixed points as a function on the manifold must have critical points.   
\medskip

  The minimal number of critical points is a topological invariant, which means that it is unchanged under homeomorphisms.  Thus, the very flexible topology of the manifold determines qualitative aspects of Hamiltonian flows. In 1983, Charles Conley and Eduard Zehnder proved this conjecture for tori of arbitrary dimension via a finite dimensional approximation of the symplectic action functional on the loop space \cite{CZ}.  The other affirmative result was due to Eliashberg in 1979, who proved it for closed two-dimensional symplectic manifolds, Riemann surfaces.   At this point, the variational methods involving finite-dimensional approximations of the action functional on the loop space stalled.  
  
Fortunately, in 1985, Mikhail Gromov pioneered the study of moduli spaces of pseudoholomorphic curves \cite{G1} to prove his celebrated nonsqueezing theorem, demonstrating that symplectic mappings are very different from volume-preserving ones. \medskip

\indent \textbf{The Gromov nonsqueezing theorem}.  A standard symplectic ball cannot be symplectically embedded into a thin cylinder. 
\medskip

Andreas Floer's subsequent breakthrough was to marry the variational methods of Conley and Zehnder with Gromov's theory of pseudoholomorphic curves, by adapting ideas from Edward Witten's interpretation of Morse theory \cite{witten}.\footnote{This brings to mind the anecdote of how Edward Witten, a physicist, came to develop his unique perspective of Morse theory.   Raoul Bott recalls first exposing Witten to Morse theory in \cite{bott}.  ``In 1979 I gave some lectures at Carg\`ese on equivariant Morse theory...to a group of very bright physicists, young and old, most of whom took a rather detached view of the lectures.  `Beautiful and oh so far from Physics' was Wilson's reaction, I remember.  On the other hand, Witten followed the lectures like a hawk, asked questions, and was clearly very interested.  I therefore thought I had done a good job indoctrinating him, so that I was rather nonplussed to receive a letter from him some eight months later, starting with the comment, `Now I finally understand Morse theory!'"

These lectures led to Witten's 1982 paper \cite{witten}, which used ideas from quantum physics to streamline Morse theory.  He recalled the evolution of these ideas in his Commemorative Lecture for the 2014 Kyoto Prize:``...Trying to get to the bottom of things, I considered simpler and simpler models, each of which turned out to contain the same puzzle. After pondering this for a long time, I eventually remembered -- I think while in a swimming pool in 1981 -- a lecture that I had heard by Raoul Bott about two years earlier...I am sure that just like me, most of the physicists at that school had never heard of it, and had no idea what it might be good for in physics. And I had probably not heard of Morse theory again until that day in 1981 when -- dimly managing to remember part of what Bott had told us -- I realized that Morse theory was behind what I had been puzzling over." }   Floer realized that the gradient trajectories counted in Morse theory didn't need to come from a flow, but instead just needed to satisfy a sufficiently nice partial differential equation with appropriate asymptotics, see \cite{F1} - \cite{F4}.  Gromov's pseudoholomorphic curves are maps between closed Riemann surfaces and symplectic manifolds that satisfy the Cauchy-Riemann equation, a nonlinear elliptic partial differential equation.  Floer modified them, studying moduli spaces of noncompact pseudoholomorphic curves perturbed by a Hamiltonian term. These \emph{Floer trajectories} are maps from the cylinder to a symplectic manifold that converge at the ends to 1-periodic solutions of the associated Hamiltonian vector field. 

At first, Gromov was skeptical of Floer's ideas.\footnote{In 1997, when Gromov was awarded the Steele Prize for a Seminal Contribution to
Research for his pseudoholomorphic curves, he recalled: ``Floer has morsified them [pseudoholomorphic curves] by breaking the symmetry, and I still cannot forgive him for this. (Alas, prejudice does not pay in science.)"  \cite{gromovfloer}. }
  Floer however, successfully formulated the nonlinear Fredholm theory describing his Floer trajectories as the zero set of an infinite dimensional bundle, thereby realizing the gradient trajectories of the highly degenerate action functional on the loop space.   This led to the creation of what is now called \emph{Floer theory}, an infinite dimensional extension of Witten's reformulation of Morse theory.    Floer used his new theory and its variants to define symplectic invariants \cite{LF} and prove the Arnold conjecture in many cases \cite{F}.

  \section{Continuing to hunt for periodic orbits with pseudoholomorphic curves}
In 1993, Hofer realized he could study moduli spaces of pseudoholomorphic maps from the complex plane to the symplectization\footnote{{The symplectization of  $(M, \mbox{ker}\alpha)$ is $(\R \times M, d(e^t \alpha))$.} } of a contact 3-manifold to prove the Weinstein conjecture for $S^3$ \cite{H1}.  However, the study of the moduli spaces of pseudoholomorphic planes is not straightforward due to additional difficulties in establishing compactness and transversality.  Clifford Taubes went on to prove the Weinstein conjecture in dimension three in 2007, relying on deep results in Seiberg-Witten theory \cite{taubes}. 

During the later '90s, Helmut Hofer, Kris Wysocki, and Eduard Zehnder continued their study of pseudoholomorphic curves in contact geometry, leading to a wealth of new dynamical results.  This work led Hofer, together with Eliashberg, to the concept of contact homology.   In 2000, constructions of these moduli based theories looked promising with the advent of a comprehensive \emph{symplectic field theory} announced in \cite{EGH}, a generalization of Floer theory and Gromov Witten theory \cite{MS2}.  This field theory involves the study of pseudoholomorphic curves from punctured Riemann surfaces to noncompact symplectic manifolds with cylindrical ends.  

These curves are still the zero set of an infinite dimensional bundle, but there is typically a failure of transversality.  As a result, one must perturb the zero set describing these curves, using either the ambient geometry or an abstract functional analytic framework.  Otherwise the resulting moduli spaces will not yield well-defined invariants.  Hofer, Wysocki, and Zehnder have developed the abstract analytic framework, collectively known as \emph{polyfolds}, to systematically resolve these issues, see \cite{HWZpoly1}-\cite{HWZgw}, and provide foundations for symplectic field theory.\footnote{This development indicates some clairvoyance on the part of George David Birkhoff, who in 1938 indicated his ``disturbing secret fear that geometry may ultimately turn out to be no more than the glittering intuitional trappings of analysis"  \cite{Bi}. On the other hand in 1980, Weinstein noted, ``The recent success of symplectic geometric methods in linear partial differential equations suggests that one might need the glitter to find the gold"  \cite{W}.}
My research, in part joint with Michael Hutchings, makes use of geometric perturbation methods  to  provide complete foundations for a subset of symplectic field theory known as \emph{cylindrical contact homology.}  These geometric methods require additional assumptions on the underlying space, but are preferable for computations and applications \cite{HN1, HN2, jo1, jo2}.

Recent work has shown that the three body problem can be studied via contact geometry \cite{P3BDY, CFK, moon}.  As a result, the modern methods of pseudoholomorphic curves are expected to give insight into the movement of satellites, allowing one to make predictions about the existence and number of energy efficient orbits that cannot be found by classical methods \cite{B}.  It would then be fitting to conclude with the words of an anonymous, albeit optimistic, symplectic geometer: ``The future of contact and symplectic geometry looks so bright that we all have to wear shades."

\end{document}